\documentclass[12pt,leqno]{amsart}

\usepackage{amsmath,amsthm,amsfonts,amssymb,latexsym,amscd,epsfig}
\usepackage{enumerate}

\usepackage{eucal}
\usepackage{euler}
\usepackage{times}

\usepackage[matrix,arrow,curve,frame]{xy}

\swapnumbers

\theoremstyle{plain}
    \newtheorem{theorem}{Theorem}[section]
    \newtheorem{lemma}[theorem]{Lemma}
    \newtheorem{corollary}[theorem]{Corollary}
    \newtheorem{proposition}[theorem]{Proposition}
    
 \newtheorem{definition-theorem}[theorem]{Definition/Theorem}

    \newtheorem*{mainthm}{Theorem}
    \newtheorem*{theorem*}{Theorem}

\theoremstyle{definition}
    \newtheorem{definition}[theorem]{Definition}

    \newtheorem{remark}[theorem]{Remark}

\theoremstyle{remark}

\numberwithin{equation}{section}

\newcounter{ictr}
\newenvironment{ilist}{\begin{list}
                         {\textup{(\alph{ictr})}}
                         {\usecounter{ictr}
                          \setlength{\leftmargin}{0.6truein}
                          \setlength{\itemsep}{0.0truein}
                          \setlength{\labelwidth}{0.3truein}}}
                      {\end{list}}

\newcommand{\R}{\mathbb{R}}
\newcommand{\C}{\mathbb{C}} 
\newcommand{\N}{\mathbb{N}}

\newcommand{\D}{\mathbb{D}}

\newcommand\op{\text{op}}
\newcommand{\ga}{\Gamma}

\newcommand{\h}{\mathcal{H}}   
\newcommand{\B}{\mathcal{B}}

\newcommand{\hyp}{\mathfrak{H}}
\newcommand{\cvx}{\mathfrak{C}}

\DeclareMathOperator{\CAT}{CAT}

\DeclareMathOperator{\st}{St}

\newcommand{\Hawaii}{Hawai\kern.05em`\kern.05em\relax i}
\newcommand{\Manoa}{M\=anoa}

\begin{document}   

\title{
Weak amenability of $\CAT(0)$-cubical groups
}  

\date{\today}

\author{Erik Guentner}
\address{University of \Hawaii, \Manoa, Department of Mathematics, 2565
  McCarthy Mall, Honolulu, HI 96822-2273}
\email{erik@math.hawaii.edu}
\author{Nigel Higson}
\address{Department of Mathematics, Pennsylvania State University,
  University Park, PA 16802}
\email{higson@math.psu.edu}

\thanks{The authors were partially supported by grants from the
U.S. National Science Foundation.}

\begin{abstract}
We prove that if $G$ is a discrete group that admits a metrically proper
action on a finite-dimensional $\CAT(0)$ cube complex $X$, then $G$ is
weakly amenable.   We do this  by constructing  uniformly bounded
Hilbert space  representations $\pi_z$  for which the quantities
$z^{\ell (g)}$ are matrix coefficients. Here $\ell$ is a length function
on $G$ obtained from the  combinatorial distance function on the complex
$X$.   
\end{abstract}

\maketitle

\section*{Introduction}
Let $G$ be a countable discrete group and denote by $C^*_\lambda (G)$
the reduced $C^*$-algebra of $G$. Every finitely supported function
$\phi\colon G\to \C$ defines a bounded linear operator $$M_\phi\colon
C_\lambda ^* (G)\to C_\lambda ^*(G)$$ by the formula 
$  (M_\phi f)(g) = \phi(g)f(g)$,
as well as bounded linear operators  
$$
M_\phi ^{(n)}\colon M_n (C^*_\lambda (G)) \to M_n (C^*_\lambda (G))
$$ 
between matrix algebras by applying $M_\phi$ to each matrix entry.  The \emph{completely bounded
  multiplier norm} of $\phi$ is defined to be the quantity  
$$ 
\| \phi \| _{\text{cb}} ={\sup} _{n} \| M^{(n)}_\phi\|.
$$
See \cite{MR996553,MR1323679}.  A countable discrete group is {\it
  weakly amenable\/} \cite{MR996553} if there exists a sequence of
finitely supported functions on $G$ converging pointwise to $1$ and
consisting of functions which are uniformly bounded in the completely
bounded multiplier norm.\footnote{This is not the standard definition,
  but see \cite{haagerup86,haagerup-kraus94}.}  Recall that a countable
discrete group $G$ is \emph{amenable} if and only if there exists a
sequence of finitely supported, positive-definite functions on $G$
converging pointwise to the constant function $1$ on $G$.  It may be
shown that the completely bounded multiplier norm of a positive-definite
function $\phi$ is equal to $\phi (e)$.  As a result, every weakly
amenable group is amenable.

The notion of amenability may be broadened in a different way. 
 A \emph{$c_0$-function} on a discrete set  is a function that extends
 to a continuous function on the one-point compactification with value
 $0$ at infinity.  A countable discrete group is said to have the
 \emph{Haagerup property} \cite{MR1852148} if there exists a   sequence of
 positive-definite $c_0$-functions on $G$ converging pointwise to the
 constant function $1$ on $G$.  Since every finitely supported function
 is a $c_0$-function, every amenable group has the Haagerup property.

Michael Cowling conjectures that every countable discrete group with the
Haagerup property is weakly amenable \cite{MR1852148}. He further conjectures
that the \emph{Cowling-Haagerup constant} of any such group is $1$,
which means that  there exists a sequence  of
finitely supported functions on $G$ converging pointwise to $1$ and
consisting of functions whose completely bounded multiplier norm is
bounded by $1$.  The conjecture is supported by the cases of
Coxeter groups  \cite{MR1925487} and  discrete subgroups of
$SO(n,1)$ \cite{decanniere-haagerup85} and $SU(n,1)$ \cite{MR748862}.   

In this paper we shall consider groups which admit an  action  on a
finite-dimensional $\CAT(0)$-cubical complex which is proper in the
sense that for every vertex $x$ in the complex, the function $d(gx,x)$
on $G$ is proper. We shall refer to  these groups as $\CAT(0)$-cubical
groups. Niblo and Reeves \cite{niblo-reeves97,MR1459140} proved that
$\CAT(0)$-cubical groups have the Haagerup property.  We shall prove the
following result: 

\begin{mainthm}
Every $CAT(0)$-cubical group is weakly amenable, with Cowling-Haagerup
constant one. 
\end{mainthm}

Our theorem adds evidence in favor of Cowling's conjecture.   In fact it
might be argued that  the class of $\CAT(0)$-cubical groups  is  quite close to
the class of all groups with the Haagerup property. It is known that  a discrete group
has the Haagerup property if 
and only if it acts properly on a {\it measured\/} space with walls
(see \cite{MR2106770} for the precise assertion).  Furthermore, a discrete
group acts on a {\it combinatorial\/} space with walls if and only if it
acts on a $\CAT(0)$ cube complex, although possibly an infinite-dimensional $\CAT(0)$ cube complex. See
\cite{chatterji-niblo04,MR2059193}.  Thus,  disregarding the issue of
infinite-dimensionality, the difference between the class of groups with
the Haagerup property and the class of $\CAT(0)$-cubical groups  comes
down to 
the difference between combinatorial and measured spaces with walls.  As
Graham Niblo 
has explained to us, the distinction is analogous to that between
trees and $\R$-trees. 

Unfortunately the possible infinite-dimensionality of $\CAT(0)$ cube
complexes can pose serious problems.   For example Farley \cite{farley}
proved that Thompson's group $F$ has the 
Haagerup property by exhibiting  a metrically proper
action of $F$ on a certain $\CAT(0)$ cube complex.  However Farley's
complex, although locally finite-dimensional, is infinite-dimensional
and our theorem above does not apply.  The question of whether or not
$F$ is   
weakly amenable remains open. 

In order to prove the weak
amenability of $\CAT(0)$-cubical groups we shall rely on a result of
Valette \cite{valette93} that makes it possible to deduce  weak
amenability from the existence of a suitable holomorphic family of
uniformly bounded representations.   Denote by  $\D$ the open unit disk
in the complex plane. 
A family $\{\pi_z\}_{z\in\D}$ of representations of a discrete group $G$
into the  bounded invertible operators on a Hilbert space $\h$ is
\emph{holomorphic}  if for every element $g\in G$ and every pair of
vectors $\xi,\xi'\in\h$ 
the matrix coefficient
\begin{equation*}
  \phi(z) = \langle\; \pi_z(g)\xi,\xi' \;\rangle 
  \end{equation*}
is a holomorphic function on $\D$.  This definition is  explored in
\cite{decanniere-haagerup85}.

\begin{theorem*}[See \cite{valette93}]
  Let $G$ be a countable discrete group.  Assume that there exists a
  proper function $\ell:\ga \to \N$ with $\ell(e)=0$ and  a holomorphic
  family $\{\pi_z\}_{z\in \D}$ of uniformly bounded representations of
  $G$ on a Hilbert space $\h$ such that $\pi_t$ is unitary for $0<t<1$.
  If there exists a vector $\xi\in \h$ such that $z^{\ell(g)} =
  \langle\, \pi_z(g)\xi,\xi \,\rangle$ for every $z\in\D$ and every
  $g\in G$, then $G$ is weakly amenable with Cowling-Haagerup constant
  one.
\end{theorem*}

We shall construct the required representations following a method of
Pimsner \cite{pimsner87} and Valette \cite{valette90} that  we  shall now
outline.  Let $X$ be a set equipped with an action of $G$ and let $\pi $
be the natural permutation representation of $G$ on $\ell^2 (X)$.  A 
function
\begin{equation*}
  c:X\times X \to \B(\ell^2 (X))
\end{equation*}
is a {\it cocycle\/} for $\pi$ if 
\begin{ilist}
  \item $c(x,x)=1$
  \item $c(v,x)c(x,y) = c(v,y)$
  \item $\pi(g) c(x,y) \pi(g)^{-1} = c(g\cdot x,g \cdot y)$
\end{ilist}
for all $v$,  $x$, $y \in X$ and all
$g\in G$.
If $c$ is a cocycle for $\pi$,  and if $x\in X$, then the formula
\begin{equation*}
  \pi_{c}(g) = c(x,gx)\pi(g)     
            \end{equation*}
defines a representation of $G$ into the  bounded invertible operators on
$\ell^2 (X)$.  If the cocycle $c$ is    
\emph{uniformly bounded}, meaning that 
$$\sup_{x,y\in X} \| c(x,y) \|<\infty,$$
then $\pi_{c}$ is uniformly bounded.  If the cocycle is \emph{unitary},
meaning that $c(x,y) = c(y,x)^*$ for all $x$ and $y$, then the
representation $\pi_c$ is unitary.

We shall take $X$ to be the set of vertices of a finite-dimensional
$\CAT(0)$ cube complex on which $G$ acts, and we shall construct a
family of cocycles $\{c_z\}_{z\in \D}$ using the geometry of $\CAT(0)$
cube complexes. In doing so we shall closely follow Pimsner and Valette,
who constructed uniformly bounded  
representations in this way when $X$ is the set of vertices of a tree
\cite{pimsner87,valette90}, 
and Januszkiewicz, who extended their results to   products of trees
\cite{januszkiewicz93}.    The same cocycle has also been analyzed very carefully by Brodzki, Niblo and Valette in unpublished work.  The construction of the cocycle is
straightforward; the proof 
that the cocycle is uniformly bounded is more difficult  and this is our  
contribution. 

\section{Overview of the Proof}
\label{overviewsec}

To prepare the reader for the more complicated case of $\CAT(0)$ cube
complexes, we shall  rapidly review the construction of Pimsner and
Valette for trees (which are simple examples of $\CAT(0)$ cube
complexes).   Let $X$ be the set of vertices of a simplicial tree on
which a group $G$ acts.  
If $z \in \D$, then define $w\in \C$ by the equation   $z^2+w^2=1$, or
$w=\sqrt{1-z^2}$, where we use the branch of the square 
root function which is holomorphic on the complement of the negative
real axis and which is positive on the positive real axis. If  
$x$, $y\in X$ are the vertices of  an edge, then define $c_z(x,y)\in
\B(\ell^2(X))$ by 
\begin{equation*}
  c_z(x,y)\delta_v = \begin{cases} 
             w\delta_x - z\delta_y, & v = x \\
             w\delta_y +z\delta_x , & v = y \\
             \delta_v, & \text{otherwise},
  \end{cases}
\end{equation*}
where $\delta_v$ denotes the Dirac function at the vertex $v$.  In
other words, $c_z(x,y)$ is given by the matrix
\begin{equation*}
  \left(\begin{matrix} w &  z \\ -z & w \end{matrix}\right) 
  \end{equation*}
on the two-dimensional subspace spanned by the ordered basis
$(\delta_x,\delta_y)$ and it is the identity on the orthogonal
complement of this subspace.  Notice that 
\begin{equation}
\label{eqn:inverse}
   c_z(x,y)=  c_z(y,x)^{-1}  ,
\end{equation}
while if $z$ is real, then $c_z(x,y)=c_z(y,x)^*$, and hence $c_z(x,y)$ is unitary.
If $x$ and $y$ are arbitrary vertices in the tree, then define
$c_z(x,y)\in \B(\ell^2 (X))$ by  
\begin{equation*}
  c_z(x,y) = c_z(v_0,v_1)  c_z(v_1,v_2) \cdots c_z(v_{n-1},v_n),
\end{equation*}
where $x=v_0,v_1,\dots,v_n=y$ are the vertices along any edge-path from
$x$ to $y$.  It follows from (\ref{eqn:inverse}) and the basic geometry
of trees that this definition is independent of the path chosen between
$x$ and $y$, and it follows from this that $c_z$ is a cocycle for $\pi$.  

Fixing a vertex $x\in X$  we obtain
representations $ \pi_z(g) = c_z(x,gx) \pi_g $.  This is a holomorphic
family --- indeed for each $g$, the operator $\pi_z(g)$ is a polynomial
in $z$ and $w$ with coefficients in $\B(\ell^2 (X))$.  Moreover, if we
define $\ell(g)$ to be $d(x,gx) $, the edge-path length between $x$ and
$gx$, then it is easy to check that $ z^{\ell(g)} = \langle
\pi_z(g)\delta_x,\delta_x\rangle $ (see for example the computations
following Lemma~\ref{techtree-lemma} below).
If $z$ is real, then $\pi_z$ is a unitary representation.  

In order to
apply Valette's theorem it remains to show that for general
$z\in \D$, the representation $\pi_z$ is uniformly bounded.  We shall
briefly indicate a series of lemmas sufficient to prove this that are
representative of what we are able to accomplish for general $\CAT(0)$
cube complexes.  Denote by $c_{ab}$ the matrix entries of $c_z(x,y)$, so
that if $b\in X$, then
$$c_z(x,y)\delta_b = \sum_{a\in X} c_{ab}\delta_a,$$
or equivalently $c_{ab} = \langle c_z(x,y)\delta_b,\delta_a\rangle$.

\begin{lemma}
\label{techtree-lemma}
  Each nonzero matrix entry $c_{ab}$ has the form  $c_{ab}=\pm z^kw^\ell$ where
$k=d(a,b)$ and where
$\ell\leq 2$.  
\end{lemma}

This is the geometric heart of the uniform boundedness argument.  For trees it may be proved by explicitly computing all the coefficients $c_{ab}$. 
Let $v_0,\dots v_n$ be the vertices on the geodesic path from $x$ to
$y$.  If  $b$ does not belong to this path, then it is easy to see that
$c_z(x,y)\delta _b = \delta _b$.  If $b=v_n$, then  
\begin{align*}
 c_z(x,y)\delta_b &= c_z(v_0,v_{1}) \cdots c_z(v_{n-1},v_{n})  \delta_{v_n} \\
    &= c_z(v_0,v_{1}) \cdots c_z(v_{n-2},v_{n-1})(w\delta_{v_{n}} + z\delta_{v_{n-1}}) \\
    &= c_z(v_0,v_{1}) \cdots c_z(v_{n-3},v_{n-2})
            (w\delta_{v_n} + wz\delta_{v_{n-1}} + z^2\delta_{v_{n-2}}) \\
    &\,\,\vdots \\
    &= w\delta_{v_n} +  wz\delta_{v_{n-1}} + wz^2 \delta_{v_{n-2}} + \cdots 
            + w z^{n-1} \delta_{v_{1}} + z^n\delta_{v_0}
\end{align*}
If $b=v_i$ with $i<n$, then  $c_z(v_{i+1},v_{i+2})\cdots c_z(v_{n-1},v_n)\delta_b=\delta_b$ and so
\begin{align*}
  c_z(x,y)\delta_b  &= c_z(v_0,v_{1}) \cdots c_z(v_{i},v_{i+1})  \delta_{v_i} \\
          &= c_z(v_0,v_{1}) \cdots c_z(v_{i-1},v_{i})(- z\delta_{v_{i+1}}+ w \delta_{v_i})\\
     &= c_z(v_0,v_{1}) \cdots c_z(v_{i-2},v_{i-1})
            (  -z\delta_{v_{i+1}} +  w^2\delta_{v_i}+zw\delta_{v_{i-1}}) \\
     &\,\,\vdots \\
          &= -z\delta_{v_{i+1}} +w^2 \delta_{v_i} + w^2  z   \delta_{v_{i-1}}
               + \cdots + w^2 z^{n-i-1}  \delta_{v_{1}} 
                    + wz^{n-i} \delta_{v_{0}},
\end{align*}
The lemma follows by inspection of these formulas.  In the general case
of cube complexes such direct computations are not so easy, and  we
shall have to work harder to obtain weaker (but adequate) results.   But proceeding with the argument for trees, we can now rapidly conclude that the cocycle $c_z$ is uniformly bounded.

\begin{lemma}
\label{tree-ub-lemma}
Let $x,y\in X$, and let $k\geq 0$.  For every   $b\in X$
$$
\#\,\left \{\, a\in X : c_{ab} = \pm z^kw^\ell \,
             \text{for some $\ell$ and all $z\in \D$}\,\right\}\le 2
$$
and for every $a\in X$
$$
\#\,\left \{\, b\in X : c_{ab} = \pm z^kw^\ell \,
           \text{for some $\ell$ and all $z\in \D$}\,\right\}\le 2.
$$ 
\end{lemma}

Both assertions follow from the fact that for any given vertex $v$ on a
geodesic interval $[x,y]$ and any $k$, there are at most two vertices on
the interval which are exactly distance $k$ away from $v$.  

\begin{proposition}
\label{tree-ub-prop}
For every $z\in\D$
  and every $x, y\in X$ the quantity  $\| c_z(x,y) \|$ is bounded by $4(1-|z|)^{-1}$.
\end{proposition}
\begin{proof}
Let $d(x,y)=N$.  Decompose
  the operator $c_z(x,y)$ on the $\ell^2$-space of the interval $[x,y]$ as a sum 
  \begin{equation*}
     c_z^{(0)} + z c_z^{(1)} + z^2 c_z^{(2)} + \dots + z^N c_z^{(N)},
  \end{equation*}
  where the operator $c_z^{(k)}$, contains only
  matrix entries of the form $\pm w^\ell$ or zero.  Since $\ell \leq 2$, each matrix coefficient of $c_z^{(k)}$
  absolute value $2$, or less.  Furthermore, by the previous lemma each
  row and column of $c_z^{(k)}$ has at most two non-zero entries.  But a
  matrix which has at most $A$ nonzero entries in each row and column,
  each of absolute value $B$ or less, determines a operator of norm $AB$
  or less.  Hence $\| c_z^{(k)} \|\leq 4$, for every $k\geq 0$, so that
\begin{equation*}
  \| c_z(x,y) \|   \le  4 \sum_0^N |z|^k,\end{equation*}
  and the proposition follows.
\end{proof}

For cube complexes we shall obtain polynomial bounds in place of the
absolute bounds in Lemma~\ref{tree-ub-lemma}, but these will be adequate
to carry through a version of the argument in the proof of
Proposition~\ref{tree-ub-prop}. 
  
In the case of trees, putting everything together, we obtain the
following result (of  Szwarc \cite{MR1092129} 
and Valette \cite{MR1316222}):  

\begin{theorem} 
  Suppose that a discrete group $G$ acts on a tree and that for some
  vertex $x$ in the tree the function $g\mapsto d(x,gx)$ is proper on
  $X$.  Then $G$ is weakly amenable with Cowling-Haagerup constant $1$.
  \qed
\end{theorem}

\section{$\CAT(0)$ cube complexes}
\label{catzerosec}

In this section we shall rapidly review the notion of cube complex and
collect the results about them that we shall need.  

A \emph{cube complex} \cite{MR919829,MR1744486,niblo-reeves98} is a set
$X$, called the 
set of \emph{vertices}, together with a collection of finite subsets of
$X$, called the \emph{cubes} of $X$, such that: every single-element set
is a cube; the intersection of any two cubes is a cube; and for every
cube $C$, there is an integer $n\ge 0$ and a bijection from $C$ to the
vertices of the standard cube $[0,1]^n$ in $\R^n$ such that the cubes in
$X$ that are subsets of $C$ correspond precisely to the sets of vertices
of the faces (of all dimensions) of the standard cube.  The cubes with
$2$ elements are called \emph{edges}; those with $4$ elements are
\emph{squares}.  A cube complex $X$ is \emph{finite-dimensional} if its
cubes are uniformly bounded in cardinality.

Every cube complex has a natural \emph{geometric realization} in which
each cube $C$ is replaced by a unit cube in $\R^n$.  This geometric
realization has a natural metric which restricts to the standard metric
on each Euclidean cube.   A  \emph{\textup{(}globally\textup{)}
  non-positively curved cube 
  complex}, or  \emph{$CAT(0)$ cube complex} is a cube complex whose
geometric realization is a $CAT(0)$ metric space \cite{MR1744486}.  For
a complex whose geometric realization is simply connected the
non-positive curvature condition can  be phrased combinatorially
\cite{MR919829,MR1744486}.   

The $CAT(0)$ condition has important consequences which we shall access
via the notion of hyperplane.    A  \emph{hyperplane} in $X$ is an
equivalence class of   edges under the equivalence relation generated by
the relation  
$$
\{x,y\}\sim \{x',y'\}\quad \Leftrightarrow \quad 
          \text{$\{x,x',y,y'\}$ is a square.}
$$
This definition is a combinatorial proxy for the notion of
\emph{geometric hyperplane} in the geometric realization of $X$: a
geometric hyperplane is the union of an equivalence class of midplanes
of cubes under the equivalence relation generated by the relation 
$$
M_1\sim M_2\quad \Leftrightarrow \quad \text{$M_1\cap M_2$ is the midplane of a cube.}
$$
See \cite{niblo-reeves98} for more information and some helpful
illustrations.  The hyperplane associated to a geometrical hyperplane is
the set of all edges bisected by the geometric hyperplane, and every
hyperplane arises in this way.  The geometric  point of view suggests
the following terminology: if an edge $\{ x,y\}$ belongs to a
hyperplane $H$, then we shall say that it  \emph{crosses} $H$; that $H$
\emph{separates} $x$ from $y$; that $x$ and $y$ are \emph{opposite} one
another across $H$;  
that $x$ and $y$ are \emph{adjacent} to $H$; and that $H$ is adjacent to
$x$ and $y$.  

The most significant consequence of the $CAT(0)$ condition is that if
$X$ is a $CAT(0)$ cube complex, then every geometric hyperplane
separates the geometric resolution of $X$ into precisely two path
components \cite[Thm~4.10]{sageev95}.  This can be expressed
combinatorially as follows.  A sequence of vertices $v_0,\dots , v_n$ in
$X$ is an \emph{edge-path} if each $\{ v_i,v_{i+1}\}$ is an edge in $X$.
The path \emph{crosses} a hyperplane $H$ if one if its edges belongs to
$H$.  Every hyperplane $H$ in a $CAT(0)$ cube complex partitions the
vertices of $X$ into precisely two sets, such that every edge path from
a vertex in one component to a vertex in the other crosses $H$, and such that
every pair of vertices in either one of the components is connected by an edge
path that does not cross $H$.  These two sets are the \emph{half-spaces}
determined by $H$ and will they be denoted $H_{\pm}$ (this involves
making an arbitrary choice of which half-space will be $H_{+}$ and which
$H_{-}$; having made such a choice we shall say that $H$ is
\emph{oriented}).  We shall say that the hyperplane $H$ \emph{separates}
two vertices from one another if one vertex lies in each of the
components $H_{\pm}$.  If $x$ and $y$ are vertices of a $CAT(0)$ cube
complex $X$, then we shall denote by $\hyp(x,y)$ the set of all
hyperplanes that separate $x$ from $y$.

The {\it distance\/} $d(x,y)$ between two vertices $x$ and $y$ of a cube
complex $X$ is the minimum length $n$ of an edge-path $v_0,\dots , v_n$
connecting them.  An edge-path from $x$ to $y$ is \emph{geodesic} if it
achieves this minimum length.

\begin{proposition}[\cite{sageev95}, Thm.~4.13]
\label{basic-crossing-prop}
  An edge-path  in a $CAT(0)$ cube complex from $x$ to $y$ crosses each hyperplane in $\hyp(x,y)$.
  An edge-path from $x$ to $y$ is   geodesic   if and only if it
  crosses only the hyperplanes in $\hyp(x,y)$, and crosses each one of
  those each exactly once. 
  \qed
\end{proposition}

It follows that $d(x,y) = \#\,\hyp(x,y)$.

\begin{corollary}
\label{cor:cvx}   
  Let $x$, $y$ and $v$ be vertices in a $CAT(0)$ cube complex.  The
  following are equivalent.  
\begin{ilist}
\item   $v$ lies on a geodesic from $x$ to $y$.
\item $d(x,y)=d(x,v) + d(v,y)$.
\item $\hyp(x,v)\cap \hyp(v,y)=\emptyset$. 
\item $\hyp(x,v)\cap \hyp(v,y)=\emptyset$ and $\hyp(x,y) = \hyp(x,v)\cup\hyp(v,y)$. 
\item $\hyp(x,v)\subseteq \hyp(x,y)$.  
\item $\hyp(x,v)\subseteq \hyp(x,y)$ and $\hyp(v,y)\subseteq \hyp(x,y)$. \qed
\end{ilist}
\end{corollary}

A {\it corner move\/} transforms an edge-path by changing a string
$u,v,w$ in the path into $u,v',w$, where $\{u,v,v',w\}$ is a square.  A
corner move does not alter the length of a edge-path or its endpoints.

\begin{proposition}[\cite{sageev95}, Thm.~4.6]
\label{prop:corner}
  Any two geodesic edge-paths in a $CAT(0)$ cube complex with the same
  endpoints differ by a sequence of corner moves. \qed
\end{proposition}

A {\it simple cancellation\/} in an edge-path replaces a
string $v,v',v$   by the singleton string $v$. 

\begin{proposition}
\label{prop:paths}
 Any two paths    in a $CAT(0)$ cube complex
  with the same endpoints are related by a sequence of corner moves and
  simple cancellations.
\end{proposition}

To prove this we shall need an additional fact about hyperplanes in
$CAT(0)$ cube complexes. Two hyperplanes $H$ and $K$ \emph{intersect} if
their geometric realizations intersect, or equivalently if there is a
square $\{x,y,z,w\}$ such that $\{x,y\}$ and $\{ z,w\}$ belong to $H$
while $\{x,z\}$ and $\{y,w\}$ belong to $K$.   

\begin{lemma}[{\cite{niblo-reeves98}, Prop.~2.10}]
\label{noselfcross-lemma} A hyperplane in a $CAT(0)$ cube complex does
not self-intersect.  In fact  if a vertex  $x$ is adjacent to a
hyperplane $H$, then there is a unique edge $\{ x,y\}$ which crosses
$H$. \qed   
\end{lemma}

We shall use the notation  $y=x^{\op}$ for the vertex $y$ opposite $x$
across $H$ (obviously this depends on the choice of $H$). 

\begin{proof}[Proof of Proposition~\ref{prop:paths}]
By the Proposition~\ref{prop:corner} it suffices to
show that any edge-path which is not a geodesic   can be reduced  to an
edge-path of shorter length  by a sequence of corner moves and simple
cancellations.   
Assume, for the sake of a contradiction, that not every non-geodesic
edge-path  is so reducible and let $v_0,\dots, v_n$ be an non-geodesic
edge-path of minimal length that  is not reducible. Since it is not a
geodesic  it must cross some hyperplane $H$ twice.  By minimality, it
must do so  between $v_0$ and $v_{1}$, and then again between $v_{n-1} $
and $v_{n}$.    Moreover the edge-path $v_1,\dots, v_{n-1}$ must be a
geodesic. Therefore by \cite[Thm.~4.13]{sageev95}, all the vertices
$v_{2},\dots, v_{n-1}$ are adjacent to $H$.  If $n=2$, then  both $v_0$
and $v_n$ are opposite $v_1$ across $H$, and so  $v_0=v_n$ by
Lemma~\ref{noselfcross-lemma}.  The path $v_0,v_1, v_2$ is therefore
$v_0,v_1,v_0$ and hence is reducible by a simple cancellation, which is
a contradiction.  If $n>2$, then by \cite[Thm.~4.12]{sageev95} the
subset  $\{v_0, v_1, v_2, v_2^{\op}\}$ is a square.  The path $v_0,
v_2^{\op},v_2, \dots , v_n$ is obtained from $v_0,\dots , v_n$ by a
corner move. But  the path  $v_2^{\op}, v_2,\dots v_n$   has length
$n-1$ and is not a geodesic (it crosses $H$ twice). It is therefore
reducible  by corner moves and simple cancellations.  It follows that
the path $v_0, v_2^{\op},v_2, \dots , v_n$  is reducible as well, and
hence so is $v_0,\dots , v_n$, which   is again a contradiction. 
\end{proof}

The \emph{convex hull} $\cvx(S)$ of a set $S$ of vertices in $X$ is the
intersection of all the half-spaces that contain $S$.  A set is
\emph{convex} if it is equal to its convex hull.  It follows from
Proposition~\ref{basic-crossing-prop} that every convex set is geodesically
convex in the sense that it contains every geodesic between any two of
its points.  In the case where $S=\{x,y\}$, the convex hull is usually
called the \emph{interval} from $x$ to $y$ and thanks to Corollary~\ref{cor:cvx} it can alternately be
characterized as the set of all points that lie on geodesics from $x$ to
$y$.

In Section~\ref{ub-sec} we shall need the following technical
result:

\begin{lemma} 
\label{cvx=interval} 
Let $x$ and $y$ be vertices in a $CAT(0)$ cube complex  and let $b$ be a vertex of some cube
$C$  that also contains $y$ as a vertex.   There is a vertex $c$ of $C$
such that $\cvx(x,y,b) \subseteq \cvx(x,c)$. 
\end{lemma}
\begin{proof} 
Let $\hyp$ be the set of all hyperplanes, necessarily passing through
the cube $C$, that separate $b$ from $y$ but not $x$ from $y$.   Let $c$
be the vertex of $C$ for which $\hyp(y,c) = \hyp$.   We shall show that   
\begin{equation*}
\cvx (x,y,b)\subseteq \cvx (x,c).
\end{equation*}
To do so, we must show that    $y\in \cvx(x,c)$ and $b\in \cvx(x,c)$.  By construction of the vertex $c$,  
$
  \hyp (x,y)$ is disjoint from $\hyp (y,c)$.  It therefore follows from Corollary~\ref{cor:cvx} that $y$ lies on a geodesic from $x$ to $c$  and hence that $y\in \cvx(x,c)$.  Similarly,   the  definition of $\hyp=\hyp(x,c)$ asserts that 
 $\hyp(y,c) = \hyp(b,c)\setminus \hyp(x,y)$,
and therefore $$\hyp(b,c)\subseteq \hyp(x,y)\cup \hyp(y,c).$$
Since $y$ lies on a geodesic from $x$ to $c$, the union is $\hyp(x,c)$. Therefore $\hyp(b,c)\subseteq \hyp(x,c)$, and so $b$ lies on a geodesic from $x$ to $c$, and hence is in the interval from $x$ to $c$.
 \end{proof}

The next three results have to do with multiple hyperplanes that are
adjacent to a single vertex.

\begin{proposition}[\cite{niblo-reeves98}, Lemma~2.14]
\label{edgemovetech-lemma}
If two hyperplanes in a  $CAT(0)$ cube complex intersect and are both adjacent to the same 
  vertex,  then they intersect in a
square containing that vertex. \qed
\end{proposition}

\begin{proposition}
\label{separatingHs-prop}
Let $H_1,\dots , H_n$ be hyperplanes in a $CAT(0)$ cube complex, all
adjacent to a vertex $x$.  If there is a vertex $y$ in $X$ such that
each hyperplane $H_i$ separates $x$ from $y$, then there is a cube of
dimension $n$ in which all the hyperplanes $H_i$ intersect.
\end{proposition}

\begin{proof}
  This follows from Proposition~3.3 in \cite{niblo-reeves98} (compare
  the remark following the proof of that proposition).
\end{proof}

\begin{proposition}
\label{prop:cube-path}
Let $x$ and $y$ be vertices of $CAT(0)$ cube complex.  There exists a geodesic path from
$x$ to $y$ that crosses all the hyperplanes in $\frak H(x,y)$ that are
adjacent to $x$ before it crosses  any other hyperplane.\end{proposition}
\begin{proof} 
  By Proposition~\ref{separatingHs-prop}, there is a cube $C$ of
  dimension $n$ containing $x$ as a vertex in which the $n$ hyperplanes
  in $\frak H(x,y)$ that are adjacent to $x$ intersect.  Let $v$ be the
  vertex diagonally across $C$ from $x$.  By the equivalence of (a) and
  (e) in Corollary~\ref{cor:cvx}, there is a geodesic from $x$ to $y$
  that passes through $v$, and any such geodesic has the required
  property. 
\end{proof}

Finally, two  hyperplanes are {\it parallel\/} if they do not intersect.

\begin{lemma}
\label{lem:half-space-intersection}
  Two parallel hyperplanes separate a $CAT(0)$ cube complex into at most three\footnote{In
  fact exactly three.} components. 
\end{lemma}

\begin{proof}
  It is convenient to work in the geometric realization.  Let $H$ and
  $K$ be parallel hyperplanes.  Since $K$ is connected, it is contained
  entirely in one component $H_{+}$ of the complement of $H$.  It
  follows that one of $K_{+}$ or $K_{-}$ does not intersect $H_{-}$, for
  if both were to intersect $H_{-}$, then there would be a geodesic
  between a point in $K_{+}$ and a point in $K_{-}$ within $H_{-}$, and
  hence a point of $K$ within $H_{-}$.  If say $K_{-}\cap
  H_{-}=\emptyset$, then
$$
\begin{aligned} X \setminus (H\cup K) &=  (H_{+}\cup H_{-})\cap (K_{+}\cup K_{-}) \\
&=  (H_+\cap K_+)\, \cup \, (H_{+}\cap K_{-}) \cup (H_-\cap K_+),\end{aligned}
$$
as required. 
\end{proof}

\section{Cocycles From Cube Complexes}

An \emph{action} of  a group on a cube complex $X$ is an action on the
set of vertices of $X$ that maps cubes to cubes. In this section we
shall associate  to the action of a group $G$ on a $CAT(0)$ cube complex
 $X$ a holomorphic family of cocycles 
$ c_z\colon X\times X \to \B (\ell^2 (X))$   parametrized   by $z\in \D$.    

As in Section~\ref{overviewsec}, we  associate to each $z\in \D$ a
complex number $w\in \C$ satisfying $z^2+w^2=1$ in such a way that $w$
depends holomorphically on $z$.   Let $x$ and $y$ be adjacent vertices
in $X$ and let $H$ be the hyperplane that separates them, oriented so
that $x\in H_{+}$ and $y\in H_{-}$.   Let us denote by $\partial H_{+}$ the set of vertices adjacent to $H$ that lie in $H_{+}$, and by $\partial H_{-}$ the set of vertices adjacent to $H$ that lie in $H_{-}$. Define a bounded operator  
$c_z(x,y)$ on $ \ell^2 (X)$ by means of the formula  
\begin{equation}
\label{protocycle-def}
  c_z(x,y)\delta_v = \begin{cases}  
        w\delta_v - z\delta_{v^{\op}}  & v\in\partial H_{+} \\
        w\delta_v + z\delta_{v^{\op}}  & v\in\partial  H_{-} \\
        \delta_v  & v\notin\partial H_{+}\cup \partial H_{-}
     \end{cases}
\end{equation}
(where $v^{\op}$ denotes the vertex adjacent to $v$ across $H$).  The definition
is consistent with the one we made for trees in
Section~\ref{overviewsec}.  However in the present, more general
situation $c_z(x,y)$ is potentially nontrivial on many more basis
vectors than $\delta_x$ and $\delta_y$ alone: it is nontrivial on the
basis vector $\delta _v$ whenever $v$ is adjacent to the hyperplane
$H$. 

Recall that a subspace of a Hilbert space is a \emph{reducing subspace}
for an operator $T$ if and only if both the subspace and its orthogonal
complement are invariant under $T$.  If the vertices $v$ and $v^{\op}$
are adjacent across $H$, then the two-dimensional subspace spanned by
$\delta _v $ and $\delta_{v^{\op}}$ is reducing for $c_z(x,y)$, as is
the subspace spanned by all $\delta_v$ for which $v$ is \emph{not}
adjacent to $v$.  Thus  by Lemma~\ref{noselfcross-lemma},  $c_z(x,y)$ is     the direct sum of a family of
operators on two-dimensional subspaces together with the identity
operator on the joint orthogonal complement of these two-dimensional
subspaces.  On the two-dimensional subspace   spanned by the ordered
pair $(\delta_v,\delta_{v^{\op}})$, where $v\in \partial H_{+}$, the
operator $c_z(x,y)$ is  given by the matrix  
\begin{equation}
\label{block}
  \left(\begin{matrix} w &  z \\  -z & w \end{matrix}\right).
\end{equation}

\begin{lemma}\label{transpose-lemma}  
  For every pair of adjacent edges
  $x$ and $y$, and for every $z\in \D$, $c_z(x,y)c_z(y,z)=I$. 
\end{lemma}
\begin{proof} 
This follows from the fact that the inverse of the matrix
(\ref{block}) is the transpose of (\ref{block}).\end{proof} 

If $x$ and $y$ are any two vertices in $X$, then we should like to define 
\begin{equation}
\label{cocycledef}
c_z(x,y) = c_z(v_0,v_1)c_z(v_1,v_2)\cdots c_z (v_{n-1},v_n),
\end{equation}
where $x=v_0, \dots , v_n = y$ are the vertices in an edge-path from $x$
to $y$.  In the case of trees, Lemma~\ref{transpose-lemma} is sufficient
to prove that the product (\ref{cocycledef}) is independent of the
choice of path from $x$ to $y$.  For general $\CAT(0)$ cube complexes we
need the following additional computation.

\begin{lemma}
\label{cornermove-lemma}
Let $x,v,y$ and $x,v',y$ be edge-paths of length two from $x$ to $y$. If the
four vertices $x,v,v',y$ span a square in the cube complex $X$, then  
$c_z(x,v)c_z(v,y)= c_z (x,v')c_z(v',y)$.
\end{lemma}
\begin{proof}
The relation between the four points is depicted in the following diagram:
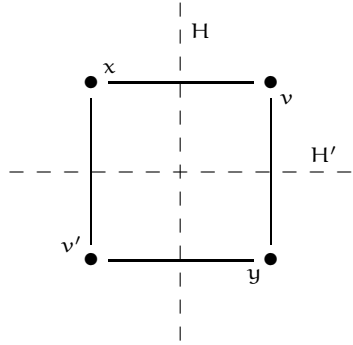
\begin{equation*}
\xymatrix{ & & \ar@{--}[dddd]^<<<<H & & \\
           & {\bullet} \ar@{-}[rr]^<x\ar@{-}[dd]_>{v'} & 
                      & {\bullet}\ar@{-}[dd]^<v & \\
                \ar@{--}[rrrr]^>>>>{H'} & & & & \\
           & {\bullet}\ar@{-}[rr]_>y & & {\bullet} & \\
                 & & & & }
\end{equation*}
The hyperplane $H$ separating $x$ and $v$ also separates $v'$ and $y$,
while the hyperplane $H'$ separating $x$ and $v'$ also separates $v$ and
$y$.  As result $c_z(x,v) = c_z(v',y)$ and $c_z(x,v')=c_z(v,y)$, since
the definition of $c_z$ for a pair of adjacent vertices depends
only on the oriented hyperplane separating them.  So our goal is
to show that
$$
c_z(x,v)c_z(x,v') = c_z(x,v')  c_z(x,v).
$$
By Lemmas~\ref{noselfcross-lemma} and \ref{edgemovetech-lemma} the  set
$X$ can be written as a disjoint union of: 
\begin{ilist}
\item Four-element sets consisting of the vertices of a square in which
  $H$ and $H'$ intersect. 
\item Two-element sets consisting of vertices adjacent across $H$ but
  with neither adjacent to $H'$. 
\item Two-element sets consisting of vertices adjacent across $H'$ but
  with neither adjacent to $H$.   
\item Single-element sets consisting of a vertex adjacent to neither $H$ nor $H'$.
\end{ilist}
The corresponding subspaces of $\ell^2 (X)$ are reducing for $c_z(x,v)$
and $c_z(x,v')$.  On the four-dimensional subspaces, the prototype of
which is the one spanned by $ \delta_x, \delta_v, \delta _{v'},
\delta_y$, the operators $c_z(x,v) $ and $ c_z (x,v')$ are represented
by matrices 
\begin{equation*}
   \left(\begin{matrix}  w & z & 0 & 0 \\
                              -z & w & 0 & 0 \\
                              0 & 0 & w & z \\
                              0 & 0 & -z & w \end{matrix}\right) 
 \quad\text{and}\quad 
                           \left(\begin{matrix}  w & 0 & z & 0 \\
                              0 & w & 0 & z \\
                             -z & 0 & w & 0 \\
                             0 & -z & 0 & w \end{matrix}\right)  
                             ,
\end{equation*}
and these two matrices  commute. On the remaining subspaces, either one
or both of $c_z(x,v)$ and $c_z(x,v')$ acts as the identity, and so the
two operators commute there too. 
\end{proof}

\begin{lemma}
  The expression (\ref{cocycledef}) defining $c(x,y)$ for general
  $x,y\in X$ is independent of the edge-path $v_0,v_1,\dots ,v_n$
  connecting $x$ to $y$. 
\end{lemma}
\begin{proof}
  According to Proposition~\ref{prop:paths} any two paths connecting $x$
  to $y$ are related by a sequence of corner moves and simple
  cancellations.  Lemma~\ref{transpose-lemma} shows that simple
  cancellations do not alter (\ref{cocycledef}).
  Lemma~\ref{cornermove-lemma} shows that corner moves do not alter
  (\ref{cocycledef}) either.
\end{proof}

\begin{proposition}
  Let $X$ be a $\CAT(0)$ cube complex and let $G$ be a
  discrete group acting on $X$.  For every
  $z\in \D$ the function $c_z:X\times X \to\B(\ell^2(X))$ defined by
  (\ref{cocycledef}) is a cocycle for the permutation representation of
  $G$ on $\ell^2 (X)$.\qed
\end{proposition} 
 
Fix $x\in X$.  Thanks to the proposition we can construct
representations $\pi_z(g) = c_z(x,gx)\pi(g)$ of the group $G$ into the bounded
invertible operators on $\ell^2 (X)$.
  
\begin{proposition}
  The family $\{ \pi_z\}_{z\in \D}$ is a holomorphic family of
  representations of $G$. If $z\in \D$ is real, then the representation
  $\pi_z$ is unitary. \qed
\end{proposition}

\begin{proposition}
  The matrix coefficient $\langle \pi_z (g) \delta _x, \delta _x\rangle$
  is equal to $z^{d(x,gx)}$.
\end{proposition}

\begin{proof}
  By definition, $\langle \pi_z (g) \delta _x, \delta _x\rangle =
  \langle c_z (x,gx) \delta _{gx} , \delta _x\rangle $, so it suffices
  to prove that
\begin{equation}
\label{matrixcoeff-fmla}
 \langle c_z (x,y) \delta _{y} , \delta _x\rangle =z^{d(x,y)}
 \end{equation}
 for every $x$ and $y$ in $X$.  Note first that by definition $c_z(x,y)$
 can be written as a product of $d(x,y)$ many operators
 $T=c_z(v_i,v_{i+1})$, each of which has the property that $\langle
 T\delta _a,\delta_b\rangle =0$ if $d(a,b)>1$.  It follows  that
\begin{equation}
\label{finiteprop-fmla} 
\langle c_z (x,y) \delta _{a} , \delta _b\rangle = 0 \quad 
       \text{if $\,\,d(a,b)> d(x,y).$}
\end{equation}
We can now prove the formula (\ref{matrixcoeff-fmla}) by induction on
$d(x,y)$.  Let $x=v_0,v_1,\dots , v_n=y$ be a geodesic edge-path from
$x$ to $y$.  Using (\ref{finiteprop-fmla}) and the induction hypothesis
we have that  
$$
c_z(v_1,v_n)\delta_{v_n} = 0\cdot \delta_{v_0} + z^{n-1}\delta_{v_{1}} + 
         \text{terms orthogonal to $\delta_{v_0},\delta_{v_1}$.}
$$
Using the explicit formula for $c_z(v_0,v_1)$ we get that
$$
c_z(v_0,v_n)\delta_{v_n} = 
     c_z(v_0,v_1)c_z(v_1,v_n)\delta_{v_n} = z^{n}\delta_{v_{0}} +
          \text{terms orthogonal to $\delta_{v_0}$},
$$
as required.
\end{proof}

\section{Calculation of Matrix Coefficients}
\label{matrix-calc-sec}

The remainder of the paper will be   devoted to proving that the
cocycle $c_z$ is  uniformly bounded.    Following the approach we took  for trees in Section~\ref{overviewsec},
we shall begin by studying  the individual matrix coefficients of the
operator $c_z(x,y)$.    Define the matrix coefficient $c_{ab}$ by the formula 
\begin{equation*}
c_z(x,y)\delta_b = \sum_{a\in X} c_{ab}\delta_a,
\end{equation*}
or equivalently $c_{ab} = \langle c_z(x,y)\delta_b,\delta_a\rangle$.  

\begin{lemma}
\label{lem:supp2}
\label{lem:supp}
\label{lem:hyp-incl}
Let $x$ and $y$ be any two vertices in $X$.   If $K$ is a hyperplane
which does not separate $x$ from $y$,  then the subspaces  $\ell^2
(K_{\pm})$    are reducing subspaces for $c_z(x,y)$.  Hence if $c_{ab}$
is non-zero for some $z\in \D$, then $\hyp(a,b)\subseteq \hyp(x,y)$. 
\end{lemma} 
\begin{proof}
Assume first that $x$ and $y$ are adjacent.  The only nonzero   matrix coefficients of $c_z(x,y)$ apart
from those on the diagonal are those $c_{ab}$ for which $a$ and $b$ are
adjacent to one another across the hyperplane $H$ separating $x$ from $y$.   So if $c_{ab}$ is nonzero for some $z\in \D$, then $a$ and $b$ must lie on opposite sides of $H$, and hence on the same side of $K$.   In general, let $v_0, \dots , v_n$ be a geodesic edge-path from $x$ to $y$, so
that  
$$
c_z(x,y)= c_z(v_0,v_1)c_z(v_1,v_2)\cdots c_z(v_{n-1},v_n) .
$$
The hyperplane $K$ separates no $v_{i}$ from $v_{i+1}$.  Therefore the subspaces $\ell^2 (K_{\pm})$ are reducing
for each $c_z(v_{i-1},v_i)$, and hence for
$c_z(x,y)$ as well. 
\end{proof}

\begin{proposition} 
\label{support-in-cvx-prop}
Let $x$, $y$ and $b$ be vertices of $X$.   If $c_{ab}$
is nonzero for some $z\in \D$, then $a$ lies in the convex hull of 
$\{ x,y,b\}$.  
\end{proposition}
\begin{proof} 
  Assume that $c_{ab}$ is nonzero for some $z\in \D$.  Let $H$ be a
  hyperplane for which all of $x$, $y$ and $b$ lie in the same
  half-space.  If $a$ lies in the other half-space then $H\in \hyp(a,b)$
  but not in $\hyp(x,y)$ contradicting Lemma~\ref{lem:hyp-incl}.  
\end{proof}

\begin{definition}
Let  $x,y\in X$.   
A {\it geodesic order\/} on $\hyp(x,y)$ is a linear order for which
there exists a geodesic edge-path   from $x$ to $y$ such that $H'<H''$
if and only if the path crosses $H'$ before it crosses $H''$.  
\end{definition}

\begin{lemma}
\label{lem:go-restrict}
If $c_{ab}$ is nonzero for some $z\in \D$, then
  every geodesic
order on $\hyp(x,y)$ induces a geodesic order on $\hyp(a,b)$.
\end{lemma}

\begin{proof}
  We shall prove the lemma by induction on the integer $d(x,y)$
  (starting with $d(x,y)=0$, where the result is trivial).   Assume that $c_{ab}\ne 0$, for some $z\in \D$.
  Let $\{ H_1, \dots, H_n\} $ be a geodesic order on $\hyp(x,y)$ and let
  $v_0,\dots, v_n$ be the corresponding geodesic path from $x$ to $y$,
  so that 
\begin{equation*}
c_z(x,y) = c_z(v_0,v_1)c_z(v_1,v_2)\cdots c_z (v_{n-1},v_n)
             =c_z(v_0,v_1)c_z(v_1,v_n).
\end{equation*}
  To prove that the given geodesic order on
$\hyp (x,y)$ restricts to a geodesic order on $\hyp (a,b)$ we shall
consider two cases.  The first is that  $a$ is not adjacent to the
hyperplane $H_1$ that separates $v_0$ from $v_1$.  In this case, since  
\begin{equation} 
\label{adj-matrix-eq}
c_{ab} = \langle c_z(v_1,v_n)\delta_b, c_z(v_0,v_1)^*\delta _a\rangle,
\end{equation} 
and since $c_z(v_0,v_1)^*\delta _a=\delta_a$,
the $ab$-matrix coefficient for $c_z(x,y)$ is equal to the $ab$-matrix
coefficient for $c_z(v_1,v_n)$ and in particular  the latter is nonzero.   By the induction hypothesis, the given
geodesic order on $\hyp(v_1,v_n)$ restricts to a geodesic order on
$\hyp(a,b)$.  But this order on $\hyp(a,b)$ is the same as the order
restricted from $\hyp(x,y)$. 
    
In the second case, $a$ is adjacent to $H_1$.  Denote by $a^{\op}$ the
vertex adjacent to $a$ across $H_1$.  From (\ref{adj-matrix-eq}) and the
definition of $c_z(v_0,v_1)$ we get that  
\begin{equation}
\label{adj-matrix-eq2}
c_{ab} = w \cdot \langle c_z(v_1,v_n)\delta_b, \delta _a\rangle \,\pm \,
z\cdot  \langle c_z(v_1,v_n)\delta_b, \delta _{a^{\op}}\rangle.
\end{equation}
 If $H_1$ separates $a$ from $b$ then by Lemma~\ref{lem:hyp-incl} the
 first inner product in (\ref{adj-matrix-eq2}) is zero, and hence 
\begin{equation*}
c_{ab} = \pm 
z\cdot  \langle c_z(v_1,v_n)\delta_b, \delta _{a^{\op}}\rangle.
\end{equation*}
Therefore the $a^{\op}b$-matrix coefficient for $c_z(v_1,v_n)$ is
nonzero. By the induction hypothesis, the given geodesic order on $\hyp
(v_1,v_n)$ restricts to a geodesic order on $\hyp(a^{\op},b)$. Since  
\begin{equation*}
\hyp(a,b) = \{ H_1\} \cup \hyp(a^{\op},b)\quad \text{and}\quad
\hyp(x,y) = \{ H_1\} \cup \hyp(v_1,v_n),
\end{equation*}
it follows easily that the geodesic order on $\hyp(x,y)$ restricts to a
geodesic order on $\hyp(a,b)$.  If $H_1$ does \emph{not} separate $a$
from $b$, then it separates $a^{\op}$ from $b$.  The second inner product in (\ref{adj-matrix-eq2}) is therefore zero, and hence \begin{equation*}
c_{ab} = 
w\cdot  \langle c_z(v_1,v_n)\delta_b, \delta _{a}\rangle,
\end{equation*}
so that  $ab$-matrix coefficient for $c_z(v_1,v_n)$ is nonzero.  By
the induction hypothesis, the geodesic order on $\hyp(v_1,v_n)$
restricts to a geodesic order on $\hyp(a,b)$.  This immediately implies
that the given  geodesic order on $\hyp(x,y)$ restricts to a geodesic
order on $\hyp(a,b)$. 
\end{proof}

\begin{lemma}
\label{lem:char-go}
A linear ordering $ \{H_1, \dots, H_n\}$ on $\hyp(x,y)$ is a geodesic
ordering if and only if the vertex $v_0=x$ is adjacent to $H_1$ and for
each $i=1,\dots ,n$ the vertex $v_i$ obtained by successively reflecting
$v_0$ across $H_1,\dots, H_{i-1}$ is adjacent to $H_i$.  In this case
the sequence of vertices $v_0,\dots , v_n$ is a geodesic edge-path from
$x$ to $y$.
\end{lemma}

\begin{proof}
If a linear ordering on $\hyp (x,y)$ is induced from a geodesic
edge-path   from $x$ to $y$, then the sequence $v_0,\dots, v_n$
is precisely the sequence of vertices along the path, and so the
adjacency condition is satisfied.  Conversely, if the adjacency
condition is satisfied, then $v_0,\dots, v_n$ is a geodesic edge-path
from $x$ to $v_n$.  The vertex $v_n$ must equal $y$ since any hyperplane
$K$ separating the two points would separate either $v_n$ from $v_0=x$
or separate $x$ from $y$.  In fact by construction of the path
$v_0,\dots, v_n$, the hyperplane $K$ would necessarily separate
\emph{both} $v_n$ and $y$ from $v_0=x$, which is a contradiction since
$K$ could not then separate $v_n$ from $y$. 
\end{proof}

\begin{lemma}
\label{lem:ells}
  Let $x,y,a,b\in X$ and suppose that no hyperplane in $\hyp(x,y)$ separates $b$ from $a$. Then  
 $\langle c_z(x,y)\delta_b,\delta_a\rangle= w ^{\ell}\langle\delta_b,\delta_a\rangle$,
where $\ell$ is the number of hyperplanes in $\hyp(x,y)$ that are  adjacent to $a$.
\end{lemma}
\begin{proof}
  The proof is by induction on $n=d(x,y)$. The case $n=0$ is trivial, so assume that $n>0$.   Let $v_0,\dots, v_n$ be a geodesic edge-path from $x$ to $y$ and let $H$ be the hyperplane that separates $v_{n-1}$ from $v_n$.   Then 
 $$
 c_z(v_{n-1},v_n)\delta_b=   \begin{cases} w\delta_b \pm z\delta_{b^{\op}}   
 & \text{if $b$ is adjacent to $H$} \\
                     \delta_{b} & \text{otherwise},
       \end{cases}
$$
where $b^{\op}$ is the vertex adjacent to $b$ across $H$.
Since $\ell^2 (H_{\pm})$ are reducing subspaces for $c_z (v_0,v_{n-1})$,  it follows that 
$ c_z (v_0,v_{n-1})\delta_{b^{\op}}$ is orthogonal to $\delta_a$.  As a result 
$$\begin{aligned}
\langle  c_z(v_{0},v_n)\delta_b,\delta_a\rangle &=  
\langle  c_z(v_{0},v_{n-1})c_z(v_{n-1},v_n)\delta_b,\delta_a\rangle \\
&= \begin{cases} w\langle c_z(v_0,v_{n-1})\delta_b,\delta_a\rangle    
 & \text{if $b$ is adjacent to $H$} \\
 \langle c_z(v_0,v_{n-1})\delta_b,\delta_a\rangle  & \text{otherwise}.
       \end{cases}\end{aligned}
$$
The result follows. \end{proof}

Suppose that $c_{ab}$ is non-zero for some $z\in \D$.  Fix a geodesic
order $\{H_1, \dots, H_n\}$ on $\hyp(x,y)$, and let $v_0,\dots,v_n$ be
the corresponding geodesic edge-path from $x$ to $y$.  Let
$$\hyp(a,b) = \{ H_{n_1},\dots , H_{n_p}\},$$
and let $a_0,\dots , a_p$
be the geodesic edge-path from $a$ to $b$, guaranteed by
Lemma~\ref{lem:go-restrict}, that corresponds to this ordering of
$\hyp(a,b)$.  Note that $d(a,b)=p$.

For $j=0,\dots, p$, let $\ell_j$ be the number of hyperplanes in 
$\{ H_k : n_j < k < n_{j+1} \}$ that are adjacent to $a_j$ (for
convenience we are setting $n_0=0$ and $n_{p+1}=n+1$).

\begin{lemma}
\label{lem:calculating-cab}
With the above notation,  $c_{ab} = \pm z^p w^{\ell_0+\dots +\ell_p}$.
\end{lemma}
\begin{proof}
We shall prove by induction on $j$, from $j=p$ down to $j=0$,  that 
\begin{equation}
\label{gen}
     \langle   c_z(v_{n_j},v_{n})\delta_b, \delta_{a_j}\rangle  =
                \pm z^{p-j} w^{\ell_j + \dots + \ell_p} . 
\end{equation}
The case $j=p$ is a consequence of Lemma~\ref{lem:ells}, while the
assertion in the current lemma is the case $j=0$.  Assume that
(\ref{gen}) holds for a given $j$.  To compute (\ref{gen}) with $j-1$ in
place of $j$ we shall write
\begin{equation*}
c_z(v_{n_{j-1}},v_{n})\delta_b = 
   c_z(v_{n_{j-1}},v_{n_j-1})c_z(v_{n_{j}-1},v_{n})\delta_b 
\end{equation*}
and then write $c_z(v_{n_{j}-1},v_{n_{p+1}})\delta_b$ as a finite linear combination
\begin{equation*}
c_z(v_{n_{j}-1},v_{n})\delta_b = \sum \alpha_k \delta_{b_k}.
\end{equation*}
As long as the sum contains no zero terms, it follows from
Lemma~\ref{lem:supp} that every vertex $b_k$ lies in the same half-space
of every $H\in \hyp (v_{n_{j-1}}, v_{n_j - 1 })$ as the vertex $b$.
Moreover $b$ and $a_{j-1}$ lie in the same half-spaces of these
hyperplanes because the only hyperplanes in $\hyp(x,y)$ that separate
$a_{j-1}$ from $b=a_p$ are $H_{n_j},\dots , H_{n_p}$.  It therefore
follows from Lemma~\ref{lem:ells} that
\begin{equation*}
\begin{aligned}
      \langle   c_z(v_{n_{j-1}},v_{n})\delta_b, \delta_{a_{j-1}}\rangle  
  & =  \sum\alpha_k \langle  c_z(v_{n_{j -1}},v_{n_j-1}) \delta_{b_k}, 
                  \delta_{a_{j-1}}\rangle \\
  &= w^{\ell_{j-1}}  \sum\alpha_k 
                \langle \delta_{b_k}, \delta_{a_{j-1}}\rangle \\
      &= w^{\ell_{j-1}}  \langle   c_z(v_{n_j-1},v_{n})\delta_b,
                                  \delta_{a_{j-1}}\rangle .
\end{aligned} 
\end{equation*}
In addition, 
\begin{equation*}
\langle   c_z(v_{n_j-1},v_{n})\delta_b, \delta_{a_{j-1}}\rangle
 =\langle   c_z(v_{n_j},v_{n})\delta_b, 
        c_z(v_{n_j-1},v_{n_j})^*\delta_{a_{j-1}}\rangle 
\end{equation*}
and 
\begin{equation*}
\langle   c_z(v_{n_j},v_{n})\delta_b,
           c_z(v_{n_j-1},v_{n_j})^*\delta_{a_{j-1}}\rangle  =
   \pm \langle   c_z(v_{n_j},v_{n})\delta_b,\bar  z\delta_{a_{j}}\rangle ,
\end{equation*}
by definition of the geodesic path $a_0,\dots , a_p$.  We conclude that 
\begin{equation*}
     \langle   c_z(v_{n_{j-1}},v_{n})\delta_b, 
                \delta_{a_{j-1}}\rangle  = 
        \pm z^{p-(j-1)} w^{\ell_{j-1} + \dots + \ell_p}
\end{equation*}
as required.
\end{proof}

\begin{lemma}
\label{symm-diff}
If $x,y,a,b\in X$, then 
$\hyp(x,y) \setminus \hyp(a,b) \subset \hyp(a,x) \triangle \hyp(b,y)$
\textup{(}symmetric difference\textup{)}.
\end{lemma}
\begin{proof}
Let $H\in \hyp(x,y) \setminus \hyp(a,b)$.  Then $x$ and $y$ are on opposite
sides of $H$ while $a$ and $b$ are on the same side.  Either $a$ and $b$ lie on 
the side containing $y$, in which  case
$H\in \hyp(a,x) \setminus \hyp(b,y)$  or they lie on the side containing
$x$, in which case  $H\in \hyp(b,y) \setminus \hyp(a,x)$.
\end{proof}

\begin{proposition}
\label{prop:cab}
  Let $x$ and $y$ be vertices of $X$.  As usual, for $a$, $b\in X$
  let $c_{ab}=\langle\, c_z(x,y)\delta_b,\delta_a \,\rangle$.  If
  $c_{ab}$ is nonzero for some $z\in \D$, then
  \begin{equation*}
    c_{ab} = \pm z^{d(a,b)} w^\ell,
  \end{equation*}
for some non-negative integer $\ell$ not exceeding twice the dimension of
$X$. 
\end{proposition}

\begin{proof} 
  In view of Lemma~\ref{lem:calculating-cab} we need only show that the
  sum $\ell_0 + \dots + \ell_p$ appearing there is bounded by twice the
  dimension of $X$.  Recall that $\ell_j$ is the cardinality of the set
  $\hyp_j$ of hyperplanes in $\{ H_k : n_j < k < n_{j+1} \}$ that are
  adjacent to $a_j$.  The only hyperplanes in $\hyp(x,y)$ that separate
  $a$ from $b$ are $H_{n_1},\dots , H_{n_p}$. Therefore
$$\hyp_0\cup \cdots  \cup \hyp_p\subseteq \hyp(x,y)\setminus \hyp(a,b),$$
and so by Lemma~\ref{symm-diff},
$$
\hyp_0\cup \cdots  \cup \hyp_p\subseteq \hyp(a,x)\, \triangle\, \hyp(b,y).
$$
We shall show that if $H$ belongs to both $\hyp_0\cup\cdots\cup \hyp_p$
and $\hyp(a,x)\setminus \hyp(b,y)$, then $H$ is adjacent to the vertex
$a$.  Similarly, we shall show that if $H$ belongs to both
$\hyp_0\cup\cdots\cup \hyp_p$ and $\hyp(b,y)\setminus \hyp(a,x)$, then
$H$ is adjacent to the vertex $b$.  It will follow from this that
$\hyp_0\cup\dots \cup \hyp_p$ is a union of two sets, the first
consisting of hyperplanes adjacent to $a$ that separate $a$ from $x$,
and the second consisting of hyperplanes adjacent to $b$ that separate
$b$ from $y$.  Since Proposition~\ref{separatingHs-prop} implies that
the hyperplanes in each set meet in a cube of $X$, it will follow that
each set can have at most $\dim(X)$ elements, and therefore that
$\ell_0+\cdots+ \ell_p$ is bounded by  $2\dim(X)$, as required.

Let $H\in \hyp_j$.  Thus $H=H_k$, where $n_j<k<n_{j+1}$, and $H$ is
adjacent to $a_j$.  Assume in addition that $H\in \hyp(a,x)\setminus
\hyp (b,y)$.  We shall show that $a$ and indeed the entire geodesic
edge-path $a_0,\dots, a_j$ is adjacent to $H$ by proving that if a
vertex $a_s$ on this path is adjacent to $H$, and if $s>1$, then
$a_{s-1}$ is adjacent to $H$ too.

Let $K$ be the hyperplane that separates $a_{s-1}$ from $a_s$ (thus
$K=H_{n_s}$.) Orient the hyperplanes $H$ and $K$ so that $b\in H_+\cap
K_+$. Then $a\in H_+\cap K_-$ since $K\in \hyp(a,b)$ while $H\notin
\hyp(a,b)$. In addition $y\in H_+$ and $x\in H_-$ since $H\notin \hyp
(b,y)$.  Since $K<H$ in the geodesic order from $x$ to $y$, there is at
least one vertex $v$ on the geodesic $v_0,\dots, v_n$ that $K$ but not
$H$ separates from $x$.  If $x\in K_-$, then $x\in H_-\cap K_-$ and
$v\in H_-\cap K_+$, while if $x\in K_+$, then $ x\in H_-\cap K_+$ and
$v\in H_-\cap K_-$.  In either case, all four of the half-space
intersections $H_\pm\cap K_\pm$ are nonempty, and it therefore follows
from Lemma~\ref{lem:half-space-intersection} that $H$ and $K$ intersect.
Since $a_{s}$ is adjacent to both $H$ and $K$, it follows from
Proposition~\ref{edgemovetech-lemma} that $H$ and $K$ intersect in a
square having $a_s$ as a vertex, and this implies that $a_{s-1}$ is
adjacent to $H$, as required.

The proof that  if $H$ belongs to both $\hyp_0\cup\cdots\cup \hyp_p$ and
$\hyp(b,y)\setminus \hyp(a,x)$, then $H$ is adjacent to the vertex $b$
is exactly the same.   
\end{proof}

\section{Uniform Boundedness of the Cocycle}
\label{ub-sec}

Our proof that the operators $c_z(x,y)$ are uniformly bounded as $x$ and $y$ range over all of
$X$ (while $z$
ranges over a compact subset of $\D$) will be based on the results of the previous section and the
following estimate of Chatterji and Ruane \cite{chatterji-ruane04} (we
shall offer our own proof of the estimate in the appendix).
 
\begin{proposition}
\label{CR-est-prop}
  Let $X$ be a finite-dimensional  $\CAT(0)$ cube complex and let $x$,
  $y$ be vertices of   $X$.  If  $k\ge 0$, then denote by $B(y,k)$ the
  set of vertices of $X$ of distance $k$ or less to $y$.  Then  
\begin{equation*}
\#\,\bigl( \cvx(x,y)\cap B(y,k)\bigr)\le (k+1)^d,
\end{equation*}
   where $d=\dim(X)$. \qed
\end{proposition}

\begin{lemma}
\label{lem:poly}
Let $x,y\in X$.  There is a polynomial function $p$ \textup{(}depending
only on the dimension of $X$\textup{)} such that for every $b\in X$ and
every $k$, 
\begin{equation*}
\#\,\left \{\, a\in X : c_{ab} = \pm z^kw^\ell  \,\text{for some $\ell$
   and all $z\in \D$}\,\right\}\le p(k)
\end{equation*}
   and such that  
  for every $a\in X$ and every $k$,
\begin{equation*}
\#\,\left \{\, b\in X : c_{ab} = \pm z^kw^\ell \,\text{for some $\ell$ and
  all $z\in \D$}\,\right\}\le p(k).
\end{equation*}
\end{lemma}

\begin{proof}
Fix $b\in X$. Let $x=u_0,\dots , u_m=b$ be a geodesic edge-path from $x$
to $b$. Let $C$ be a cube containing the vertex $b$ in which meet all
the hyperplanes that are adjacent to $b$ and that separate  $b$ from
$y$. Let  $b=v_0,\dots, v_n=y$  be a geodesic edge-path from $b$ to $y$
composed of a path $v_0,\dots, v_i$ consisting of vertices in $C$,
followed by a path $v_{i+1},\dots v_n$   which does not cross any
hyperplane adjacent to $b=v_0$ (see Proposition~\ref{prop:cube-path}).
By the cocycle property,  
\begin{equation*}
   c_z(x,y)= c_z(u_0,u_1)\cdots c_z(u_{m-1},u_m)c_z(v_0,v_1)\cdots c_z(v_{n-1}, v_n).
\end{equation*}
 Furthermore $c_z(v_j,v_{j+1})\delta_b=\delta _b$ for all $j\ge i$ since $b$ is not adjacent to the hyperplane separating $v_j$ from $v_{j+1}$, and therefore
$c_z(x,y)\delta_b = c_z(x,v_i)\delta_b$. Now, because $b$ and $v_i$ are
vertices of the same cube $C$, it follows from
Proposition~\ref{support-in-cvx-prop} and Lemma~\ref{cvx=interval} that
there is another vertex $c$ of $C$ such that  
\begin{equation*}
\left\{\, a\in X : c_{ab}\ne 0  \,\text{for some $z\in \D$}\,\right\}
\subseteq \cvx(x,c).
\end{equation*}
Note next that $d(b,c)\le \dim(X)$, so that  if $d(a,b) =k$, then
$d(a,c)\le k+\dim(X)$.  In addition, if  $c_{ab} = \pm z^kw^\ell$, then we
proved in Proposition~\ref{prop:cab} that $k=d(a,b)$, from which it
follows that $d(a,c)\le k+\dim(X)$. 
Therefore 
\begin{equation*} 
\left \{\, a\in X : c_{ab} = \pm z^kw^\ell  \,\text{for some $\ell$ and all
    $z\in \D$}\,\right\}\subseteq \cvx(x,c)\cap B(c,k+\dim(X))
\end{equation*}
and Proposition~\ref{CR-est-prop} implies that for all $k$,
\begin{equation*}
 \#\, \left \{\, a\in X : c_{ab} = \pm z^kw^\ell  \,\text{for some $\ell$ and
 all $z\in \D$}\,\right\}\le  (k+d+1)^d, 
\end{equation*}
where $d=\dim(X)$.  
 
To prove the second estimate, it suffices to note that 
$c_z(x,y) = c_{\bar z} (y,x)^*$,
so that   the $ab$-matrix coefficient of $c_z(x,y)$   is equal to   the
complex conjugate of the $ba$-matrix coefficient of $c_{\bar z} (y,x)$.
Thus   the second estimate, using the same polynomial $p(k) =
(k+d+1)^d$, follows from the first. 
\end{proof}

Following the approach we took in Section~\ref{overviewsec}, let us use
the results of Section~\ref{matrix-calc-sec} to decompose   the operator
$c_z(x,y)$ as a linear combination 
\begin{equation*}
\label{eq:ksum}
  c_z(x,y) = \sum_{k\ge 0} z^k  c^{(k)}_z(x,y),  
\end{equation*}
in which  the non-zero matrix coefficients of the operators
$c^{(k)}_z(x,y)$ are all of the form  $\pm w^\ell$, where $0\le \ell\le 2
\dim(X)$.  Note that this is actually a finite linear combination since,
by Lemma~\ref{lem:hyp-incl},
$c_z^{(k)}(x,y)=0$ when $k\ge d(x,y)$.

\begin{proposition}
For every compact subset $K\subseteq \D$, 
\begin{equation*}
  \sup\,\bigl\{\,  \|c_z(x,y)\|\, :\, x,y\in X,\, z\in K\,\bigr \}<\infty.
\end{equation*}
\end{proposition}
\begin{proof} 
  Since $|w|^2\le 2$, the matrix entries of $ c_z^{(k)}(x,y)$ are all
  bounded in absolute value by $2^{ \dim(X)}$.  Furthermore there are
  at most $p(k)$ 
  non-zero matrix entries in each row and column, where $p$ is the
  polynomial function of Lemma~\ref{lem:poly}.  Therefore, as in
  Section~\ref{overviewsec}, 
\begin{equation*} 
\| c_z^{(k)}(x,y)\| \le 2^{ \dim(X)}p(k) 
\end{equation*}
It follows that 
\begin{equation*}
\|c_z(x,y)\| \le  2^{ \dim(X)} \sum_{k=0}^\infty |z|^k p(k) 
\end{equation*}
and this gives the result. 
\end{proof}

Having established that the cocycle $c_z(x,y)$ is uniformly bounded, the
proof of our main theorem is complete:

\begin{definition}  
  A \emph{$CAT(0)$-cubical group} is a group $G$ which admits an
  action
  on a finite-dimensional $CAT(0)$ cube complex
  in such a way that $d(gx,x)$ is a proper function on $G$ for some (and
  hence any) vertex $x$.
\end{definition}

\begin{theorem}
If $G$ is a  $CAT(0)$-cubical group, then $G$  is weakly amenable and has
  Cowling-Haagerup constant $1$. \qed
\end{theorem}

\appendix

\section{A Bound on the Intersections of Intervals with Balls}

We shall give a new proof, which may be of independent interest,  of the following result of Chatterji and Ruane
\cite{chatterji-ruane04}.   

\begin{theorem}
\label{thm:RD-est}
  Let $X$ be a $\CAT(0)$ cube complex and let $x$, $y$ be vertices of
  $X$.  For every $r\in\N$ the cardinality of the set $\cvx(x,y)\cap
  B(x,r)$ is bounded by $(r+1)^d$, where $d$ is the dimension of $X$.
\end{theorem}

The proof relies on the
following proposition.

\begin{proposition}
\label{prop:dec}
Let $X$ be a $\CAT(0)$ cube complex and let $x$ and $y$ be vertices of
$X$.  Let $H\in\hyp(x,y)$ and suppose that $x$ is adjacent to $H$.
There exists a second vertex $v$ adjacent to $H$ and on the same side of
$H$ as $x$ for which
\begin{equation*}
\label{decomp}
  \cvx(x,y) = \cvx(x,v)\cup \cvx(x^{\text{\rm op}},y) 
           \quad\text{\textup{(}disjoint union\textup{)}}.
\end{equation*}
\end{proposition}


\begin{remark}
  In fact we only require the forward inclusion for Theorem~\ref{thm:RD-est}.
\end{remark}

We  shall also use the fact that  a hyperplane $H$ in a $CAT(0)$ cube complex 
may be given a natural $\CAT(0)$ cube complex structure in its own
right. The set underlying $H$ can be taken to be the set of vertices in
$\partial  H_{+}$, for any fixed orientation of $H$.  The cubes in $H$
are exactly the subsets of $H$ that are cubes in $X$.    Distance   and
convex hulls  may be computed in $H$ or in the ambient complex $X$.  See
\cite{sageev95}.   

\begin{proof}[Proof of Theorem~\ref{thm:RD-est}, assuming Proposition~\ref{prop:dec}]
  The proof is by double induction on the dimension $d$ of $X$ and the
  radius $r$.  The base of the induction comprises two cases.  The
  result is obvious in the case of arbitrary $d$ and $r=0$ since
  $B(x,0)$ contains only $x$.  It is also obvious in the case of
  arbitrary $r$ and $d=0$. 

For the induction step, given $\bold d>0$ and $\bold r>0$,  assume the
  result for all complexes of dimension less than $\bold d$, no matter
  what the value of $r$, and for all complexes $X$ of dimension $d=\bold
  d$ and all balls in $X$ of radius $r<\bold r$.  Now assume that
  $\dim(X)=\bold d$ and let $r=\bold r$.  Let $x,y\in X$ and let $H$ be
  a hyperplane adjacent to $x$ that separates $x$ from $y$.  It 
  follows from the proposition that 
\begin{equation*}
  \cvx(x,y) \cap B(x,r) 
       \subseteq \bigl( \cvx(x^{\op},y) \cap B(x^{\op},r-1) \bigr)
                     \cup \bigl( \cvx(x,v) \cap B(x,r) \bigr) 
\end{equation*}
(note that if $u\in \cvx(x^{\op},y)$ then $d(x^{\op},u) = d(x,u)-1$
since $H$ separates $x$ from $u$ but not $x^{\op}$ from $u$).  The
induction hypothesis implies that
\begin{equation*}
  \# \; \bigl( \cvx(x^{\op},y) \cap B(x^{\op},r-1)\bigr)  \leq r^d.
\end{equation*}
Since $\cvx(x,v)\subseteq H$ and since the dimension of $H$ is less than
$\bold d$, the induction hypothesis also implies that
\begin{equation*}
  \#\;\bigl ( \cvx(x,v) \cap B(x,r)\bigr ) \leq (r+1)^{d-1}.
\end{equation*}
Combining these estimates, we conclude \begin{equation*}
  \#\; \left ( \cvx(x,y) \cap B(x,r) \right ) \leq r^d + (r+1)^{d-1} \leq (r+1)^d,
\end{equation*}
as required.  \end{proof}

The remainder of the appendix is devoted to proving
Proposition~\ref{prop:dec}.  Let $x$ and $y$ be vertices of $X$ and let
$H\in \hyp(x,y)$ be such that $x\in \partial H$.  Orient $H$ by the
requirement that $x\in H_+$.  Let $v$ be an element of the finite set
$\cvx(x,y)\cap H_+$ of maximal distance from $x$ (or equivalently of
minimal distance to $y$).  We shall show that $v$ has the properties
stated in the proposition.

\begin{lemma}
\label{lem:thin}
Every vertex in  $\cvx(x,y)\cap H_+$ is adjacent to $H$.  In
particular, $v$ is adjacent to $H$.
%
%
\end{lemma}
\begin{proof}
  Let $u\in \cvx(x,y)\cap H_+$ and suppose for the sake of a
  contradiction that $u$ is not adjacent to $H$.  Let $u_0,\dots , u_n$
  be a geodesic from $x$ to $y$ containing $u$.  There is a first vertex
  $u_j$ on this geodesic that is not adjacent to $H$.  Then of course
  $u_{j-1}$ is adjacent to $H$. Let $K$ be the hyperplane separating
  $u_{j-1}$ from $u_j$. Both $H$ and $K$ separate $u_{j-1}$ from $y$.
  It follows that there is a square containing $u_{j-1}$ as a vertex in
  which $H$ and $K$ intersect, and as in the proof of
  Proposition~\ref{prop:cab}, it follows that $u_j$ is adjacent to $H$.
  Contradiction. \end{proof}

\begin{lemma}
If $K\in \hyp(v,y)$ and if $K$ is adjacent to $v$, then $K=H$.  
\end{lemma}
\begin{proof}
  If $K\neq H$, then since both $H$ and $K$ separate $v$ from $y$, and
  since both are adjacent to $v$, the two hyperplanes intersect in a
  square containing $v$ as a vertex.  The vertex $v^{\op}$ adjacent to
  $v$ across $K$ is therefore adjacent to $H$.  It is also in the
  interval from $x$ to $y$ and further away (by one) from $x$ than $v$.
  This contradicts the definition of $v$.
\end{proof}

In the next lemma we require a small amount of the theory of {\it normal
  cube paths\/}.  
See \cite{niblo-reeves98} for further information, especially the remark
following Proposition~3.3.

\begin{lemma}
\label{lem:disjoint}
  No hyperplane $K\neq H$ separating $v$ and $y$ can intersect $H$.
\end{lemma}
\begin{proof}
Assume for the sake of a contradiction that such a hyperplane $K$ exists. 
There is then a hyperplane  $K$ other than $H$ that
separates $v$ from $y$ and which has the property that if $K$ intersects
the normal cube $C$, then no hyperplane (other than $H$ itself)
intersecting a normal cube prior to $C$ on the normal cube
path from $v$ to $y$ intersects $H$.  Observe
that the cube $C$ cannot be the first cube in the normal cube path from $v$
to $y$ since by the previous lemma $K$ cannot be adjacent to $v$.   Let
$K_1,\dots,K_d$ be the hyperplanes 
  spanning the normal cube $C'$ 
  immediately preceeding $C$.    A simple separation argument based on
  Lemma~\ref{lem:half-space-intersection} shows that $K\cap K_i$ is
  nonempty for each $i$.  Thus, the hyperplanes $K,K_1,\dots,K_d$
  intersect pairwise, and each is adjacent to the (unique) vertex $w$ in
  $C'\cap C$.  By Proposition~\ref{edgemovetech-lemma} each pair
  intersects in a square having $w$ as a vertex.  But, the link of $w$
  is a flag complex (see \cite[Thm.~II.5.20]{MR1744486}), so that
  these hyperplanes
  intersect in a cube having $w$ as a vertex.  This contradicts the
  definition of normal cube path according to which
$\st(C')\cap C= \{\, w \,\}$.
\end{proof}

\begin{proof}[Proof of   Proposition~\ref{prop:dec}]
We shall show that 
\begin{equation}
\label{eq:minus}
      \cvx(x,y)\cap H_-  = \cvx(x^{\op},y)
\end{equation}
and 
\begin{equation}
\label{eq:plus}
      \cvx(x,y)\cap H_+  = \cvx(x,v).
\end{equation}
Let $u\in \cvx(x,y)\cap H_-$.  Then $d(x,u)+d(u,y) = d(x,y)$. Since
  $u\in H_-$ it follows that   
  $d(x^{\op},u)=d(x,u)-1$.  Therefore
  \begin{equation*}
    d(x^{\op},u)+d(u,y) = d(x,u) -1 + d(u,y) = d(x,y)-1 = d(x^{\op},y)
  \end{equation*}
so that $u\in \cvx(x^{\op},y)$.  The other inclusion in (\ref{eq:minus})
is obvious.

It follows immediately from \cite[Thm.~4.13]{sageev95} that
$\cvx(x,v)\subseteq H_+$.  Using Corollary~\ref{cor:cvx}, if
$w\in\cvx(x,v)$ then $\hyp(x,w)\subseteq\hyp(x,v)\subseteq\hyp(x,y)$ so
that $w\in \cvx(x,y)$.

Finally, let  $w\in \cvx(x,y)\cap H_+$. Then by Corollary~\ref{cor:cvx},
\begin{equation*}
  \hyp(x,w)\subseteq \hyp(x,y) = \hyp(x,v)\cup\hyp(v,y) 
            \quad \text{(disjoint union)}.
\end{equation*}
We want to show that  
$\hyp(x,w)\subseteq \hyp(x,v)$, or equivalently that
$\hyp(x,w)\cap\hyp(v,y)$ is empty.  This is a separation argument.  
Indeed, suppose $K\in \hyp(v,y)$.  If $K=H$ then $K\notin \hyp(x,w)$ by
assumption.  If $K\neq H$ then 
$K$ and $H$ are parallel by Lemma~\ref{lem:disjoint}.  Further, $K$ is
contained entirely in $H_-$, since $K\in\hyp(v^{\op},y)$ and both
$v^{\op}$ and $y\in H_-$.  Now, a
geodesic path from $x$ to $w$ is completely contained in $H_+$, so
cannot cross $H$, and therefore cannot cross $K$.  That is,
$K\notin\hyp(x,w)$.  
\end{proof}


\bibliographystyle{amsalpha}
\bibliography{cube_references}

\end{document}